\documentclass[12pt]{article}
\usepackage{amssymb}

\usepackage{amsmath,amsthm,bbm}
\usepackage{version}
\usepackage{color}

\author{
John E. M\raise.45ex\hbox{c}Carthy
\thanks{Partially supported by National Science Foundation Grant  
DMS 2054199. These are the personal views of the author, and do not represent the opinion of the National Science Foundation.}
}

\title{Conferences---an Owner's Manual}

\begin{document}
\maketitle

I have just attended my first in-person conferences since the start of COVID. They form a critical part of our profession, so I decided to write my opinions about how to maximize their benefits. I say nothing of online conferences, because I don't understand them.

\section{Personal History}

The first out-of-town conference I attended was at the University of Arkansas in April 1988.
There was a mini-course of 5 lectures by Harold Shapiro on Quadrature Domains, and individual lectures by 
other senior operator theorists and complex analysts (in those days, graduate students rarely traveled to conferences, and never gave talks). My adviser, Donald Sarason, had arranged for several of his students to go, and to share a room in a hotel
(a Hilton! I had heard of the luxury hotel brand, but had never set foot in one, let alone slept in one.)

This conference turned out to be the most important event of my professional life. I found it terribly exciting---real mathematicians talking about their work and progress on interesting problems. I could understand the statements, even if the proofs were complicated. Outside the talks, the professors talked to the graduate students as if we were real people.
Allen Shields even invited me to join a group for dinner. (Unfortunately I was too shy to accept---a foolish mistake).
At the conference banquet, I sat at a table with several professors who chatted with us, and Philip Davis  gave a fairly long after-dinner speech that was witty, interesting, and for me the high-light of the evening. 

Before I went to the conference I was feeling somewhat desperate. I was in my fourth year as a Ph.D. student, had no results, and every time I sat down to work the sense that I had to prove something RIGHT NOW caused my anxiety and blood pressure to spike, making progress even more difficult. I was starting to believe that I would never finish my degree.

After I came back, I was newly enthused about mathematics, and within two weeks had proved a theorem that turned into my thesis. 
After my adviser signed off on my thesis, with the psychological burden lifted, I became a much better mathematician\footnote{It is bizarre that you spend 5 years as a Ph.D. student, essentially writing your first paper, and then you are told that to be competitive you now have to write
3 papers a year. Amazingly you can. If  someone could only write down whatever the magical knowledge is
that one acquires writing a dissertation, Ph.D. completion times could be reduced to a year.}. 
Quadrature domains, which I had never heard of before the mini-course, have turned up in my professional life in unexpected and interesting ways. Many of the people I met at that conference became friends. 

Since then, I have attended numerous conferences and organized several. Some have been more enjoyable than others.
Here are my personal opinions on what makes for a good conference.

\section{The participants}

Mathematics is a human 
activity (see \cite{jxm} for my thoughts on this), and a social activity.
We may prove theorems on our own, but we need to communicate with others. Conferences serve several social 
functions, the relative importance of which change over one's career.
\begin{itemize}
\item They serve to educate---we learn from the talks,  we learn from discussions, and we learn from conversations at meals. We learn both mathematics, and mathematical culture.
\item They serve to advertise---here is my new theorem! Let me explain why it is interesting!
\item They serve to socialize. When you go to your first conference, you probably only know a few people from your own university. But over time you get to keep meeting people you have met before, some of whom become friends, even close friends.
\item They serve to network. This is like socializing, but there is a subtext of helping you professionally. After all, you are much more likely to get a job offer from a university if someone on the faculty there has seen you talk.
\item They serve to inspire. I am really impressed by the theorem the speaker is telling us about. Next year I want to be on the stage talking about my own impressive theorem!
\item
They serve to enthuse. Mathematics is really interesting and fun!
\end{itemize}

Conference attendees are a heterogeneous group---graduate students, postdocs, junior faculty, senior faculty, some mathematicians who are close to the core theme of the conference, some who are quite far from it, occasionally undergraduates and non-professionals. They all bring different things, and want different things from the conference. For the conference to be successful, they must all cooperate.

\section{Conference Talks}

There are several excellent articles on giving mathematical talks; see for example the essays by  B. Kra \cite{kra12} and T. Tao \cite{ta09} on talks in general, and the post by W. Ross \cite{ros08} on 20 minute talks in particular. My views on colloquium talks are here \cite{mcc99a}. In this essay, I will confine myself specifically to conference talks.

\begin{itemize}
\item
The best medium for most mathematics talks is a blackboard, perhaps with some interruptions for graphics.
This is for two reasons. The mathematical one is that blackboards allow far more material to be visible, so the audience
can check back on definitions and statements. The psychological one is that it forces the speaker to proceed slowly, since it takes time to write. (See V. Peller's essay \cite{pe12} for the advantages of the blackboard).

\item
At many conferences blackboards are unavailable. If you are giving a Beamer talk:

\begin{itemize}
\item GO SLOWLY and don't show a lot of writing.

\item
Theorem: 
For an $N$ minute talk,  the optimal number of slides\footnote{You get the title page for free.} is $\frac N2$.

Any number larger than $N$ is malpractice.
It takes much longer for the audience to absorb ideas and statements than most speakers realize\footnote{I find my
ability to absorb decreases over the course of the conference, as I get more tired.}.

\item
Never put a full paragraph on a slide. Write the minimum necessary---it does not  have to be in full sentences.
We are all conditioned to read whatever is put in front of us\footnote{This is why billboards are a driving hazard.}. Time spent reading is time not spent listening to the speaker.
\item
Do not end with a slide that just says ``Thank You''. Your last slide should contain your main result(s), so that audience members can look at it while absorbing your talk and thinking about questions.\footnote{You can put ``Thank You" on the last line if you want---a good use of the {\em \\pause} command in Beamer.}

\end{itemize}

\item
There may well be people in the audience who heard you talk on a similar theme before. Don't let this influence your presentation. Many other people in the audience haven't attended your talk before, and even the ones who did attend
don't remember very much about it.

\item
Know the mathematical range of the audience. Try to make the talk worthwhile for all of them, not just a couple of experts.

\item
Everybody in the audience has chosen to attend your talk, instead of spending their time proving a theorem (or whatever else
humans do when they are not attending math talks). It is your obligation to make sure  their faith in you is justified, and the time at your talk is well-spent. 

\item
A good talk takes a lot of preparation. Don't cheat your audience by not preparing properly. Do not use slides or notes from another talk---make a fresh preparation for this particular audience. Of course, there may be significant overlaps with prior talks, but your emphasis should shift depending on who the audience is (and what the theme of the conference is).
Two thirds of a good 60 minute talk is not a good 40 minute talk.

\item
Don't go over time. The chair of the session should give the speaker a sign when there is 5 minutes remaining, and stand up when the time is up\footnote{Boris Korenblum famously unplugged the projector when a speaker went over time and wouldn't stop.}.

\item
Iff you like the talk, tell the speaker. Everybody likes positive feedback. And if nobody is coming up to you after your talk to tell you they liked it, perhaps you should wonder what you should be doing differently.

\end{itemize}

\section{Social Behavior}

\begin{itemize}
\item
Go to talks, not just by well-known mathematicians. 

\item
Most conferences have plenary lectures (this means no other talk is scheduled simultaneously) and parallel sessions. 
The purpose of plenary talks is to inform---they should be like colloquia, but aimed at the audience designated by the conference\footnote{So at an AMS meeting, a plenary talk should be just like a colloquium. At a conference on Hilbert Function spaces, the speaker can assume that the audience already knows what Hilbert Function spaces are, and believes them to be inherently interesting. But the speaker still has 
to convince the audience that their particular subtopic is  interesting.}.
Plenary talks close to you educate you about your current area of research. Plenary talks far from your current interests may
educate you about your future area of research, or at least help you see a bigger picture of where your work fits.

\item
Parallel sessions have a range of speakers, from senior mathematicians to people giving their first ever conference talk.
Attend talks in the parallel sessions too.
Remember, everybody needs an audience. You want people to come to your talk, don't you?

\item
Sadly, not every speaker will have taken the lessons of the previous section to heart, and some talks will be boring.
If you get lost in a talk, it is perfectly acceptable to take out a pad of paper and work on your own mathematics.
But don't type on a keyboard---this is distracting for the audience around you.

\item
Questions at the end of talks are great, but it takes time to formulate a good question, especially if the talk covered a lot of material.
Sometimes the chair of the session will tap somebody in advance and ask them to ask the first question\footnote{
Sometimes it is hard to ask a question because the speaker anticipated the natural questions and answered them during the talk.
Here is a good question if it hasn't already been answered: What is a very simple application of your theorem?}.

\item
Talk to people, at coffee breaks and meals. Go out to eat with other attendees at lunch and dinner.\footnote{Do this in manageable groups.
Having 20 odd mathematicians walk down the street looking for a restaurant for lunch where they can all sit together is ridiculous. Once the group exceeds 6 or 8 you can only talk to those sitting close to you anyway, so you may as well start in a smaller group.}

\item
Introduce yourself to people you don't already know, and talk to them. This is hard, especially if you are junior, but do it anyway. Some of these people will become your friends and collaborators, some will tell you interesting stories. 

\item
If it is a themed conference, there will probably be a banquet. Go to the banquet. Try and ensure that some of the people
sitting at your table are people you don't already know\footnote{If you already know everybody at the conference, good for you! I have nothing to teach you.}.

\item
Conference Organizers: Try to have the banquet in a venue where everybody at a table can hear everybody else near them.
This means not too much background noise, and either small round tables (8 people maximum, 6 is better) or rectangular tables.

\item
It is traditional for the first speaker at the conference to give a speech at the banquet\footnote{
If for some reason they cannot give the speech, they should arrange for somebody else to do so.}.
This speech doesn't have to be long, but it should have some message in addition to thanking the organizers.

\item
Organizing a conference is a lot of work. Don't make unreasonable requests of the organizers.

\item
Organizing a conference is a lot of work. Thank the organizers! 

\end{itemize}

\section{Summary}

Conferences are central to the practice of mathematics. Their success relies on your efforts as an attendee. This includes being a good audience member, being a good
 speaker if you have the privilege of speaking at the conference, and making an effort to socialize, especially with people you don't know, or who seem to need some help gaining entree. Done right, conferences are enjoyable and stimulating, and you will go home with a renewed enthusiasm to do mathematics.

\bibliographystyle{plain}

\bibliography{../references_uniform_partial}
\end{document}